\documentclass[12pt]{article}
\usepackage{amsfonts}

\newtheorem{thm}{Theorem}[section]

\newtheorem{cor}[thm]{Corollary}

\newtheorem{prop}[thm]{Proposition}
\newtheorem{conjecture}{Conjecture}

\newenvironment{remark}{\par\medskip\noindent{\bf Remark.\ }}{\par\smallskip}

\newcommand{\be}{\begin{equation}}
\newcommand{\ee}{\end{equation}}
\newcommand{\openbox}{\leavevmode
  \hbox to8pt{\hfil\vrule\vbox to6pt{\hrule width6pt\vfil\hrule}\vrule}}

\newcommand{\qed}{\hbox to5pt{ } \hfill \openbox\bigskip\medskip}

\newcommand{\rk}{\mbox{\rm rank}}

\newcommand{\ve}[1]{\mathbf{#1}}

\newcommand{\cF}{\mbox{$\cal F$}}

\newcommand{\cH}{\mbox{$\cal H$}}

\newcommand{\cV}{\mbox{$\cal V$}}

\newcommand{\R}{\mathbb R}

\title{A new upper bound for codes with a single Hamming distance}
\author{G\'abor Heged\"{u}s
\\{\normalsize  \'Obuda University}
\\{\normalsize B\'ecsi \'ut 96/B, Budapest, Hungary, H-1032}
\\{\normalsize hegedus.gabor@uni-obuda.hu}
}

\begin{document}
\maketitle

\begin{abstract}
In this short note we give a new upper bound for the size of a set family with a single Hamming distance.

Our proof is an application of the linear algebra bound method. 
\end{abstract}
\medskip

\section{Introduction}

Throughout the paper 
$n$ denotes a positive integer and $[n]$ stands for the set $\{1,2,
\ldots, n\}$. We denote the family of all subsets of $[n]$ by $2^{[n]}$. 

Let $\cF\subseteq 2^{[n]}$ be a  family of subsets and $F, G\in \cF$ be two distinct elements of $\cF$.
Let $d_H(F,G)$ denote the Hamming distance of the sets
 $F$ and  $G$, i.e.,  
$d_H(F,G):=|F \Delta G|$, 
where $F \Delta G$ is the usual symmetric difference.

Let $q>1$ be an integer. Let $\cH\subseteq \{0,1,\ldots ,q-1\}^n$ and 
let $\ve h_1,\ve h_2\in \cH$ be two elements of the vector system $\cH$.  
Let $d_H(\ve h_1,\ve h_2)$ stand for the Hamming distance of the
vectors $\ve h_1,\ve h_2\in \cH$:
$$
d_H(\ve h_1,\ve h_2):=|\{i\in [n]:~  (\ve h_1)_i\neq (\ve h_2)_i\}|.
$$

Delsarte proved the following remarkable upper bound for the size of the vectors systems with $s$ distinct Hamming distances in \cite{D},  \cite{D2}. Babai, Snevily and Wilson gave a different proof of this upper bound, which used the polynomial method.

\begin{thm} \label{Delsarte}
Let $0<s\leq n$, $q>1$ be positive integers. Let $L=\{\ell_1 ,\ldots ,\ell_s\}\subseteq [n]$ be a set of $s$ positive integers. Let $\cH\subseteq \{0,1,\ldots ,q-1\}^n$ and suppose that $d_H(\ve h_1,\ve h_2)\in L$ for each distinct $\ve h_1,\ve h_2\in \cH$ vectors. Then
$$
|\cH|\leq \sum_{i=0}^s {n\choose i} (q-1)^i.
$$
\end{thm}

In the $q=2$ special case we get the following statement.

\begin{cor} \label{Delsarte2}
Let $0<s\leq n$ be positive integers.  Let $L=\{\ell_1 ,\ldots ,\ell_s\}\subseteq [n]$ be a set of $s$ positive integers. 
Let $\cF\subseteq 2^{[n]}$ be a set system such that $d_H(F,G)\in L$ for each distinct $F,G\in \cF$. Then
$$
|\cF|\leq \sum_{j=0}^s {n\choose j}.
$$
\end{cor}

Our main result is inspired by the nonuniform version of  Fisher's inequality (see \cite{M}, \cite{I}).  To the best of the author's knowledge, this result is completely new.

\begin{thm} \label{main}
Let $\cF=\{F_1, \ldots , F_m\}$ be  a family of subsets of $[n]$ such that there exists a $\lambda>0$ positive integer with $d_H(F_i,F_j)=\lambda$ for each $i\neq j$. Suppose that $\lambda\neq \frac{n+1}{2}$. Then $m=|\cF|\leq n$.
\end{thm}
%\proof
It is easy to verify that our upper bound is  sharp. We describe here only one important example.

Let $n=4k$ and if it exists, then let $H$ denote a Hadamard matrix of order $4k$.  Denote by ${\ve h}(i)$  the $i^{th}$ row of the matrix $H$. Let ${\ve k}(i)\in \{0,1\}^n$  denote  the corresponding (0-1) vector, i.e., ${\ve k}(i)_t=0$ if and only if  ${\ve h}(i)_t=-1$. 

Consider  the set system $\cF$, which corresponds to the set of characteristic vectors $\{ {\ve k}(i):~ 1\leq i\leq n\}$. Then clearly $d_H(F,G)=\frac n2$ for each distinct $F,G\in \cF$ and $|\cF|=n$ (see e.g. \cite{J} Theorem 15.7).

%Let $q$ be a prime power and let $\cP:=PG(2,q)$ denote the projective plane over $\Fq$. In this example $\cP$ is our base set.  Clearly $|\cP|=q^2+q+1$. Let $\cF$ denote the set system of projective lines of $\cP$. Then $|\cF|=|\cP|$ and $d_H(F,G)=2q$ for each distinct $F,G\in \cF$.
%We conjecture that  and $s\neq n-1$. 

In Section 2 we prove  first Theorem \ref{main}, then we give a compact proof for the special case of  Theorem \ref{main}, when $\lambda\leq \frac{n}{2}$.
In Section 3 we present a conjecture, which gives a possible generalization of our main result for $q>2$.

\section{Proofs}

The proof of  our main inequality is based on  the linear algebra bound method and the following result on the determinant of a special matrix.

\begin{prop} \label{det}
Let $J_m$ be the $m\times m$ all one matrix and let $I_m$ denote the  $m\times m$ identity  matrix. Let $\theta,\gamma\in \R$ be fixed real numbers. Then
$$
\det(\theta J_m+\gamma I_m)=(\gamma +m\theta )\gamma^{m-1}.
$$
\end{prop}
Proposition \ref{det} was proved in \cite{BF} Exercise 4.1.3.\\

{\bf Proof of Theorem \ref{main}:}\\

Let $M$ denote the $m\times n$ $(-1,1)$ {\em signed incidence matrix} of the set system $\cF$, i.e., define $M(i,j)$ as
\[ M(i,j) =\left\{ \begin{array}{ll}
1 & \textrm{if $j\in F_i$} \\
-1 & \textrm{otherwise.}
\end{array} \right. \]

It is easy to verify that we can summarize the  condition $d_H(F_i,F_j)=\lambda$ for each $i\neq j$ in the matrix equation 
$$
N:=M\cdot M^T=(n-2\lambda)J_m+2\lambda I_m.
$$

It follows from Proposition \ref{det} with the choices $\theta:=n-2\lambda$, $\gamma:=2\lambda$  that
$$
\det(N)=\det((n-2\lambda)J_m+2\lambda I_m)=(2\lambda+m(n-2\lambda))(2\lambda)^{m-1}.
$$

Hence $\det(N)\neq 0$ precisely if $2\lambda+m(n-2\lambda)\neq 0$. 
Clearly if $\det(N)\neq 0$, then 
$$
m=\rk(M\cdot M^T)\leq \rk(M)\leq n,
$$
and we proved our result.

On the other hand it follows from Corollary \ref{Delsarte2} that $m\leq n+1$. 
Our aim is to prove that $m\leq n$. Assume the contrary and  suppose that $m=n+1$. Then $\det(N)=0$, hence it follows from  $2\lambda+m(n-2\lambda)=0$ that $\lambda=\frac{n+1}{2}$, a contradiction. \qed

\begin{remark}
If $\lambda\leq \frac{n}{2}$, then it is easy to verify that the matrix $N$ is positive definite.
\end{remark}

{\bf Proof of Theorem \ref{main} in the special case $\lambda\leq \frac{n}{2}$:}

Let $\ve v_i$ denote the $(-1,1)$ {\em signed} characteristic vector of the set $F_i$, i.e, 
\[ (\ve v_i)_j =\left\{ \begin{array}{ll}
1 & \textrm{if $j\in F_i$} \\
-1 & \textrm{otherwise.}
\end{array} \right. \]

It is enough to prove that $\ve v_1, \ldots ,\ve v_m\in \R^n$ are linearly independent vectors over $\R$.
 
We prove an indirect way. Assume that there exists such  linear relation 
$$
\sum_{i=1}^m \mu_i \ve v_i=\ve 0,
$$
where not all coefficients $\mu_i$ are zero. 

It follows from  the  condition $d_H(F_i,F_j)=\lambda$ for each $i\neq j$ that
$\langle \ve v_i,\ve v_j\rangle=n-2\lambda$ if $i\neq j$ and clearly $\langle \ve v_i,\ve v_i\rangle=n$ for each $i$.

Consider the expression
$$
A:=\langle \sum_{i=1}^m \mu_i \ve v_i,\sum_{j=1}^m \mu_j \ve v_j\rangle =0.
$$
Clearly then
$$
A=\sum_{i=1}^m  \mu_i^2  \langle \ve v_i,\ve v_i\rangle+ \sum_{1\leq i\ne j\leq m} \mu_i\mu_j \langle \ve v_i,\ve v_j\rangle\geq  
$$
$$
\geq \sum_{i=1}^m  \mu_i^2 \cdot n +\sum_{1\leq i\ne j\leq m} \mu_i\mu_j \cdot (n-2\lambda).
$$
But 
$$
\sum_{i=1}^m  \mu_i^2 \cdot n +\sum_{1\leq i\ne j\leq m} \mu_i\mu_j \cdot (n-2\lambda)=2\lambda\sum_{i=1}^m  \mu_i^2 +\sum_{i=1}^m  \mu_i^2 \cdot (n-2\lambda)+ \sum_{1\leq i\ne j\leq m} \mu_i\mu_j \cdot (n-2\lambda)=
$$
$$
=2\lambda\sum_{i=1}^m  \mu_i^2 + (n-2\lambda)\cdot (\sum_{i=1}^m  \mu_i)^2
$$
and 
$$
2\lambda\sum_{i=1}^m  \mu_i^2 +(n-2\lambda) \cdot (\sum_{i=1}^m  \mu_i)^2>0,
$$
because $\lambda>0$, $n-2\lambda\geq 0$ and there exists an index $i$ such that the coefficient $\mu_i$ is nonzero by the indirect condition. Consequently $A>0$, a contradiction. \qed

\section{Concluding remarks}

We conjecture the following generalization of Theorem \ref{main}.
\begin{conjecture} \label{Hconj}
Let $\cV=\{\ve v_1, \ldots , \ve v_m\}\subseteq \{0,1,\ldots ,q-1\}^n$ be  a vector system   such that there exists a $\lambda>0$ positive integer with $d_H(\ve v_i,\ve v_j)=\lambda$ for each $i\neq j$. Suppose that $\lambda\neq \frac{(q-1)(n+1)}{q}$. Then $m=|\cV|\leq n(q-1)$.
\end{conjecture}

It would be very interesting to characterize the extremal configurations   in Theorem \ref{main}, i.e. the set systems $\cF=\{F_1, \ldots , F_m\}$  and the $\lambda>0$ positive integers such that $d_H(F_i,F_j)=\lambda$ for each $i\neq j$ and  $|\cF|=n$.
%{\bf Acknowledgements.} 
Finally it would be desirable to determine if Theorem \ref{main}  remains valid in the case $\lambda=\frac{n+1}{2}$.

\end{document}